\documentclass[12pt,reqno]{amsart}
\usepackage{amssymb,amscd,amsbsy}
\usepackage{amsthm}
\newcommand{\lb}{\linebreak}

\newcommand{\nl}{\newline}

\renewcommand{\a}{\alpha}

\newcommand{\g}{\gamma}

\newcommand{\vk}{\varkappa}
\newcommand{\z}{\zeta}

\newcommand{\s}{\sigma}

\newcommand{\f}{\varphi}

\newcommand{\w}{\omega}

\renewcommand{\L}{\Lambda}

\renewcommand{\O}{\Omega}
\newcommand{\U}{\Upsilon}

\newcommand{\A}{{\mathcal A}}

\newcommand{\E}{{\mathcal E}}

\newcommand{\cL}{{\mathcal L}}
\newcommand{\M}{{\mathcal M}}

\newcommand{\V}{{\mathcal V}}
\newcommand{\W}{{\mathcal W}}

\newcommand{\1}{{\bf 1}}

\newcommand{\C}{{\Bbb C}}
\newcommand{\T}{{\Bbb T}}
\newcommand{\pp}{{\Bbb P}}
\newcommand{\dd}{{\Bbb D}}

\newcommand{\Z}{{\Bbb Z}}
\newcommand{\mm}{{\Bbb M}}
\newcommand{\0}{{\boldsymbol{0}}}

\newcommand{\bs}{\boldsymbol}

\newcommand{\rf}[1]{(\ref{#1})}

\newcommand{\df}{\stackrel{\mathrm{def}}{=}}
\newcommand{\dist}{\operatorname{dist}}
\newcommand{\Ker}{\operatorname{Ker}}

\newcommand{\spn}{\operatorname{span}}

\newcommand{\rank}{\operatorname{rank}}

\newcommand{\eeq}{\end{equation}}
\newcommand{\beq}{\begin{equation}}
\newcommand{\bay}{\begin{eqnarray}}
\newcommand{\ey}{\end{eqnarray}}
\newcommand{\bey}{\begin{eqnarray*}}
\newcommand{\eey}{\end{eqnarray*}}

\newcommand{\be}{\infty}

\newcommand{\bl}{\blacksquare}
\newcommand{\ess}{\operatorname{ess}}
\newcommand{\ind}{\operatorname{ind}}
\newcommand{\wind}{\operatorname{wind}}
\newcommand{\Range}{\operatorname{Range}}

\newcommand{\Pf}{{\bf Proof. }}

\newcommand{\ov}{\overline}

\newtheorem{thm}{\hspace{\parindent}Theorem}[section]

\newtheorem{lem}[thm]{\hspace{\parindent}Lemma}

\theoremstyle{remark}

\newtheorem*{rem*}{Remark}

\newcommand\cE{\mathcal{E}}

\begin{document}

\newcommand{\vse}{\vspace{.2in}}
\numberwithin{equation}{section}

\title{\bf Very badly approximable matrix functions}
\author{V.V. Peller and S.R. Treil}
\thanks{The first author is partially supported by NSF grant DMS 0200712.
The second author is partially supported by NSF grant DMS 0200584 }
\maketitle

\begin{abstract}
We study in this paper very badly approximable matrix functions on 
the unit circle $\T$, i.e., matrix functions $\Phi$
such that the zero function is a superoptimal approximation of $\Phi$. 
The purpose of this paper is to obtain 
a characterization of the continuous very badly approximable functions.

 Our characterization is more geometric than
algebraic characterizations earlier obtained in 
\cite{PY} and \cite{AP}. It involves analyticity of certain families of subspaces defined in terms
of  Schmidt vectors of the matrices $\Phi(\z)$, $\z\in\T$. 
This characterization can be extended to the wider class of {\em
admissible} functions, i.e., the class of matrix functions
$\Phi$ such that the essential norm
$\|H_\Phi\|_{\rm e}$ of the Hankel operator $H_\Phi$
is less than the smallest nonzero superoptimal 
singular value of $\Phi$.

In the final section we obtain a similar
characterization of badly approximable matrix functions. 
\end{abstract}

\setcounter{section}{0}
\section{\bf Introduction}
\setcounter{equation}{0}

\

A well-known classical result in complex analysis says that any bounded measurable function $\f$ on
the  unit circle $\T$ has a best uniform approximation by bounded analytic functions, i.e., there exists a
function $f\in H^\infty$ such that 
$$
\|\f - f\|_\infty = \dist_{L^\infty}(\f, H^\infty)= \inf_{h\in H^\infty} \|\f-
h\|_\infty. 
$$
It is even more remarkable, that in many cases
the best approximation $f$ is unique. For example, this is true if $\f$ is continuous on $\T$;
this was first proved for the first time in \cite{Kh} and was rediscovered later by several other
mathematicians.

A function $\f\in L^\infty$ is called {\em badly approximable} if 
$$
\|\f\|_\infty = \dist_\infty(\f, H^\infty),
$$
i.e., if its norm cannot be reduced by subtracting an $H^\infty$ function.
Another way to describe badly approximable functions is to say that any such
function is the difference between a function and its best approximation in
$H^\infty$.

There is an elegant characterization of the set of continuous 
badly approximable functions:
a nonzero continuous function $\f\not\equiv0$ on the unit circle $\T$
is badly approximable if and only if it has constant modulus and
its winding number $\wind\f$ is negative (see \cite{AAK}, \cite {Po}). 
Recall that the winding number of a continuous function $\f:\T \to
\C\setminus\{0\}$,
 is  the number of turns of the point $\f(e^{it})$ around the origin when $t$ runs
from $0$ to $2\pi$ (see, e.g., \cite{P}, Ch. 3, \S 3).

This characterization can be extended to broader classes of functions, for which the
winding number is not defined. For such functions the result can be stated in terms of
Hankel and Toeplitz operators. 

\newcommand{\coker}{\operatorname{coker}}

It is well known (see e.g., \cite{D}) that if $\f\in C(\T)$ and $\f$ does not vanish on $\T$, \lb then 
the Toeplitz operator $T_\f$ on the Hardy class $H^2$ is Fredholm and \lb$\ind T_\f
=-\wind\f$ (recall that for a Fredholm operator $A$, its {\it index} is defined as
$\dim\Ker A -\dim \Ker A^*$). The above characterization of badly approximable
functions can be easily generalized in the following way: if
$\f$ is a function in
$L^\be$ such that the essential norm
$\|H_\f\|_{\rm e}$ of the Hankel operator $H_\f$ is less than its norm, then $\f$
is badly approximable if and only if $\f$ has constant modulus almost  everywhere
on $\T$, $T_\f$ is Fredholm, and $\ind T_\f>0$ (see e.g., \cite{P}, Ch. 7, \S 5). 
Recall that the {\it Toeplitz operator}
$T_\f:H^2\to H^2$ and the {\it Hankel operator} $H_\f:H^2\to H^2_-\df L^2\ominus H^2$ are defined by
\bay
\label{th}
T_\f f=\pp_+\f f,\quad H_\f f=\pp_-\f f,
\ey
where $\pp_-$ and $\pp_+$ are the orthogonal projections onto the subspaces $H^2$ and $H^2_-$.
Recall also that 
$$
\|H_\f\|=\dist_{L^\be}(\f,H^\be)\quad\mbox{and}\quad\|H_\f\|_{\rm e}=\dist_{L^\be}(\f,H^\be+C)
$$
(see, e.g. \cite{P}).

\medskip

{\bf 1.1. Badly approximable matrix functions.} In this paper we deal with matrix-valued
functions. The notion of a badly approximable matrix function can be defined in a similar way. A
matrix function
$\Phi$ with values in the space
$\mm_{m,n}$ of $m\times n$  matrices is called {\it badly approximable} if 
$$
\|\Phi\|_{L^\be}=\inf\{\|\Phi-F\|_{L^\be}:~F\in H^\be(\mm_{m,n})\}.
$$
Here 
$$
\|\Phi\|_{L^\be}\df\ess\sup_{\z\in\T}\|\Phi(\z)\|_{\mm_{m,n}},
$$
$\mm_{m,n}$ is equipped with the standard operator norm, and $H^\be(\mm_{m,n})$ is the space of
bounded analytic functions with values in $\mm_{m,n}$.

While it is possible (and it is done in this paper) to describe  badly
approximable matrix-functions, the problem does not look very natural. The main
reason is, that even for continuous matrix-valued functions a best  $L^\infty$ approximation by
analytic matrix functions is almost never unique. For example,
suppose that $m=n=2$ and suppose that $u$ is a scalar 
badly approximable unimodular function (i.e., $|u(\z)|=1$ almost everywhere on $\T$). 
Consider now the matrix function
$\Phi=\left(\begin{array}{cc}u&\bs{0}\\\bs{0}&\0\end{array}\right)$. it is easy to see that for any 
scalar function $f$ in the unit ball of $H^\be$, the matrix function 
$\left(\begin{array}{cc}\0&\bs{0}\\\bs{0}&f\end{array}\right)$ is a best approximation of $\Phi$.
Clearly, if  $\psi$ is an arbitrary 
scalar function in the unit ball of $L^\be$, then the matrix function
$\left(\begin{array}{cc}u&\bs{0}\\\bs{0}&\psi\end{array}\right)$
is badly approximable. However, $\psi$ can as ``bad'' as possible.

The problem of describing all badly approximable functions such that $\0$ is the
unique best approximation looks slightly more natural. This problem is also solved in
this paper, see Theorem \ref{t5.2} below. But the most natural problem
appears when one considers the approximation method that gives a unique
``very best'' approximation (for continuous matrix-valued functions).

Thus in our opinion, in the case of matrix functions it is most natural to consider
the notion of very badly approximable matrix functions, which was introduced in
\cite{PY}. To define a very badly approximable  matrix function, we need the notion
of superoptimal approximation (see \cite{PY}). 

\medskip

{\bf 1.2. Superoptimal approximations and very badly approximable matrix functions.}
Recall that for a matrix $A$ the {\it singular value} $s_j(A)$, $j\ge0$, is, by 
definition, the distance from $A$  to the set of matrices of rank at
most $j$. Clearly, $s_0(A)=\|A\|$.

\medskip

{\bf Definition.}
Given a matrix function $\Phi\in L^\be(\mm_{m,n})$ we define inductively
the sets $\bs{\O}_j$, $0\le j\le\min\{m,n\}-1$, by
$$
\bs{\O}_0=\{F\in H^\be(\mm_{m,n})
:~F~\mbox{minimizes}~\ t_0\df\ess\sup_{\z\in\T}\|\Phi(\z)-F(\z)\|\};
$$
$$
\bs{\O}_j=\{F\in \O_{j-1}:~F~\mbox{minimizes}~\ 
t_j\df\ess\sup_{\z\in\T}s_j(\Phi(\z)-F(\z))\},\quad j>0.
$$
Functions in $\bigcap\limits_{k\ge 0} \bs{\O}_k = \bs{\O}_{\min\{m,n\}-1} $ are called {\it
superoptimal approximations} of $\Phi$ by bounded analytic matrix functions. The numbers
$t_j=t_j(\Phi)$ are called the {\it superoptimal singular values} of $\Phi$. Note that the functions
in $\bs{\O}_0$ are just the best approximations by analytic matrix functions.
%(a function $F\in H^\be(\mm_{m,n})$ is called a {\it best approximation} of $\Phi$ if
%$\|\Phi-F\|_{L^\be}=\inf_{G\in H^\be(\mm_{m,n})}\|\Phi-G\|_{L^\be}$).

A matrix function $\Phi$ is called {\it very badly approximable} if the zero function is 
a superoptimal approximation  of $\Phi$. Again, a very badly approximable function  can be
interpreted as the difference between a function and its superoptimal approximation.

\medskip

{\bf 1.3. Some known results.}
 The notion of superoptimal approximation seems very natural for the approximation theory of
matrix-valued functions, for the superoptimal approximation is unique for continuous functions:
it was designed to have uniqueness!  Namely, it was shown In
\cite{PY} that if
$\Phi\in(H^\be+C)(\mm_{m,n})$ (i.e., all entries of
$\Phi$ belong to $H^\be+C$), then $\Phi$ has a unique superoptimal approximation $Q$ by bounded
analytic matrix  functions. 
Moreover, it was shown in
\cite{PY} that
\bay
\label{sv}
s_j(\Phi(\z)-Q(\z))=t_j(\Phi)\quad\mbox{for almost all}\quad\z\in\T.
\ey

The problem to describe the very badly approximable functions was posed in \cite{PY}. It follows from
\rf{sv} that if $\Phi$ is a very badly approximable function in \lb$(H^\be+C)(\mm_{m,n})$, then the singular values 
$s_j(\Phi(\z))$ are constant for almost all $\z\in\T$. Moreover, it was shown in \cite{PY} that if 
in addition to this $m\le n$ and
\mbox{$s_{m-1}(\Phi(\z))\ne0$} almost everywhere, 
then the Toeplitz operator \lb\mbox{$T_{z\Phi}:H^2(\C^n)\to H^2(\C^m)$}
has dense range (if $\Phi$ is a scalar function, the last condition is equivalent to the fact that 
$\ind T_\Phi>0$).
Note that the Toeplitz and the Hankel operators whose symbols are matrix functions can be 
defined in the same way as in the scalar case (see \rf{th}). Obviously, this necessary condition is equivalent to the 
condition $\Ker T_{\bar z\Phi^*}=\{\0\}$. In fact, the proof of necessity given in \cite{PY} allows one to obtain 
a more general result: if $\Phi$ is an arbitrary very badly approximable function in $(H^\be+C)(\mm_{m,n})$ and
$f\in\Ker T_{\bar z\Phi^*}$, then $\Phi^*f=\0$.

On the other hand, in \cite{PY} an example of a continuous $2\times2$ function $\Phi$ was given
such that
$s_0(\Phi(\z))=1$,
$s_1(\Phi(\z))=\a<1$, $\z\in\T$, $T_{z\Phi}$ is invertible but $\Phi$ is not even badly approximable.

The very badly approximable matrix functions of class $(H^\be+C)(\mm_{m,n})$ were characterized in
\cite{PY}  algebraically, in terms of so-called thematic factorizations. 

Later in \cite{PT} the above results of \cite{PY} were generalized to 
the broader context of matrix functions
$\Phi$ such that the essential norm $\|H_\Phi\|_{\rm e}$ of the Hankel 
operator $H_\Phi$ is less than the smallest
nonzero superoptimal singular value of $\Phi$. We call such matrix functions $\Phi$ {\it admissible}.
In particular, if $\Phi$ is an admissible very badly approximable $m\times n$ matrix function, then
the functions $s_j(\Phi(z))$ are constant almost everywhere on $\T$ and
$$
\Ker T_{\bar z\Phi^*}=\{f\in H^2(\C^n):~\Phi^*f=\0\}.
$$

In \cite{AP} another algebraic characterization of the set of very badly approximable 
admissible matrix functions was given in terms of canonical factorizations (see \S \ref{s2}
for the definition).

We refer the reader to the book \cite{P}, which contains all the above information and 
results on superoptimal approximation and very badly approximable functions.

\medskip

{\bf 1.4. What is done in the paper.}
Although a complete description (necessary and sufficient condition) of very badly
approximable matrix functions was obtained in \cite{PY} and \cite{AP}, this description is rather
complicated: it says that a function is very badly approximable if and only if it admits some
special factorization. While such characterizations are very helpful for {\em constructing} very
badly approximable functions, it is not easy to check, using such characterizations, that a function
is very badly approximable.

The main result of the paper is Theorem \ref{plo} in \S \ref{s4}, which gives another 
 description of admissible   very badly  approximable matrix-functions.  In particular,
it gives a complete description of the very badly approximable matrix-functions with entries in
$H^\infty+C$. This description is more  geometric and closer 
in spirit to the scalar result stated at the beginning of this paper
than the algebraic characterizations obtained in \cite{PY} or \cite{AP}. 

Note, that the result is new and highly nontrivial even for continuous
functions. The main difficulty is to understand  the structure of very
badly approximable functions, not to extend the results to a wider class
of functions.

The paper is organized as follows: 
In \S \ref{s3} we find a new necessary condition for an admissible matrix functions to be
very badly approximable. It involves analyticity of certain families of subspaces. 
However, we will see in 
\S \ref{s3} that 
if we add this analyticity condition to the above two necessary conditions, 
the three conditions will still remain
insufficient.

In \S \ref{s4} we slightly modify this necessary conditions to obtain a description of
the very badly approximable admissible matrix functions. In \S 5 we give a new approach
to this problem that is based on the notion of a superoptimal weight.

Finally, in \S \ref{s5} we obtain a characterization of the set badly approximable matrix
functions $\Phi$ satisfying the condition $\|H_\Phi\|_{\rm e}<\|\Phi\|_{L^\be}$ and we
obtain a characterization of badly approximable matrix functions, for which $\0$ is the unique best
approximation.

In \S \ref{s2}  we define canonical factorizations and state several results we are going to
use in \S \ref{s3} and later to establish the main result of the paper.

\medskip

{\bf1.5. Acknowledgement.} The first author is grateful to I.M. Gelfand and \lb M. Atiyah for 
encouraging conversations.

\medskip

{\bf 1.6. Notation.} Throughout this paper we use the following notation:
\nl
$I_n$ is the identity matrix of size $n\times n$;
\nl
$\bs{I}_n$ is the matrix function on $\T$ equal to $I_n$ almost everywhere;
\nl
$\bs{0}$ denotes a scalar or matrix function on $\T$ that is equal to zero almost everywhere;
\nl
$\bs{1}$ is the scalar function identically equal to $1$.
\nl
$z$ denotes the identical function: $z(\z)=\z$, $\z\in\T$.

\

\setcounter{section}{1}
\section{\bf Preliminaries}
\label{s2}
\setcounter{equation}{0}

\

To define canonical factorizations, we need the notion 
of balanced unitary-valued functions. Recall that a matrix function $G\in H^\be(\mm_{m,n})$ is called {\it inner} if 
\lb$G^*G=\bs{I}_n$. A matrix function $G\in H^\be(\mm_{m,n})$ is called {\it outer} if $GH^2(\C^n)$ is dense in 
$H^2(\C^m)$. Finally, $G\in H^\be(\mm_{m,n})$ is called {\it co-outer} if the transposed function
$G^{\rm t}\in H^\be(\mm_{n,m})$ is outer.

It is easy to deduce from the definition of co-outer functions that if $G$ is a co-outer function in  $H^\be(\mm_{m,n})$
and $f$ is a function in $L^2(\C^n)$ such that $Gf\in H^2(\C^m)$, then $f\in H^2(\C^n)$ (see e.g., \cite{P}, Ch. 14, \S 1).

By the Beurling--Lax--Halmos theorem (see e.g., \cite{N}), a nonzero subspace $\cL$ of $H^2(\C^n)$ is
invariant under multiplication by $z$ if and only if $\cL=\U H^2(\C^r)$, where $1\le r\le n$ and $\U$ is an inner 
$n\times r$ matrix function. It is easy to see that
\bay
\label{r}
r=\dim\{f(\z):~f\in\cL\}\quad \mbox{for almost all}\quad\z\quad\mbox{in the unit disk}\quad\dd.
\ey

\medskip

{\bf Definition.} Let $n$ be a positive integer and let $r$ be an integer
such that $0<r<n$. Suppose that $\U$ is an $n\times r$ inner and co-outer matrix 
function and $\Theta$ is an $n\times(n-r)$ inner and co-outer matrix function. If
the matrix function
$$
\V=\left(\begin{array}{cc}\U&\ov{\Theta}\end{array}\right)
$$
is unitary-valued, it is called an {\it$r$-balanced matrix function}. 
If $r=0$ or $r=n$, it is natural 
to say that an $r$-balanced matrix is a constant unitary matrix.
An $n\times n$ matrix function $\V$ is called {\it balanced} if it is $r$-balanced 
for some $r$, $0\le r\le n$. $1$-balanced matrix functions are also called {\it thematic}.

\medskip

It is well known (see \cite{V}) that each inner and co-outer matrix function $\U$ has a balanced completion 
$\left(\begin{array}{cc}\U&\ov{\Theta}\end{array}\right)$.

The following result was obtained in \cite{AP}.

\medskip

{\bf Theorem A.} {\it Let $\V$ be a balanced matrix function. Then the Toeplitz operators
$T_\V$ and $T_{\V^{\rm t}}$ have trivial kernel and dense range.}

\medskip

We also need the following fact from \cite{AP}.

\medskip

{\bf Theorem B.} {\it Suppose that $\Phi\in L^\be(\mm_{m,n})$ and 
{\em$\|H_\Phi\|_{\text e}<\|H_\Phi\|$}. Let $\cL$ be the minimal invariant
subspace of multiplication by $z$ on $H^2(\C^n)$ that contains all maximizing 
vectors of $H_\Phi$. Then
$$
\cL=\U H^2(\C^r),
$$
where $r$ is the number of superoptimal singular values of $\Phi$ equal to 
$\|H_\Phi\|$ and $\U$ is an inner and co-outer $n\times r$ matrix function.}

\medskip

If we apply Theorem A to the transposed function $\Phi^{\rm t}$, we find an $m\times r$ inner and
co-outer matrix function $\O$ such that the invariant subspace of multiplication by $z$ on $H^2(\C^m)$
spanned by all maximizing vectors of $H_{\Phi^{\rm t}}$ coincides with $\O H^2(\C^r)$.

Consider now balanced completions $\left(\begin{array}{cc}\U&\ov{\Theta}\end{array}\right)$
and $\left(\begin{array}{cc}\O&\ov{\Xi}\end{array}\right)$ of $\U$ and $\O$ and define the unitary-valued functions
$\V$ and $\W$ by
$$
\V=\left(\begin{array}{cc}\U&\ov{\Theta}\end{array}\right),\quad
\W=\left(\begin{array}{cc}\O&\ov{\Xi}\end{array}\right)^{\rm t}.
$$

\medskip

{\bf Theorem C.} {\it Under the hypotheses of Theorem A the matrix functions $\U$, $\Theta$, $\O$, $\Xi$
are left invertible in $H^\be$.}

\medskip

Recall that a matrix function $\Phi\in H^\be(\mm_{m,n})$ is said to be {\it left invertible in $H^\be$} 
if there exists $\Psi\in H^\be(\mm_{n,m})$ such that $\Psi\Phi=\bs{I}_n$. Theorem C was established in [AP], see also 
\cite{PT} where it was proved in the case when $\V$ and $\W^{\rm t}$ are $1$-balanced.

The following result can also be found in \cite{AP}.

\medskip

{\bf Theorem D.} {\it Suppose that $\Phi\in L^\be(\mm_{m,n})$ and $\|H_\Phi\|_{\rm e}<\|H_\Phi\|$. 
Then $\Phi$ is badly approximable if and only if it admits a factorization
$$
\Phi=\W^*\left(\begin{array}{cc}\s U&\bs{0}\\\bs{0}&\Psi\end{array}\right)\V^*
$$
where $\V$ and $\W^{\rm t}$ are $r$-balanced matrix functions, $r$ is the number of superoptimal
singular values of $\Phi$ equal to $\|\Phi\|_{L^\be}$, $\s=\|\Phi\|_{L^\be}$,
$U$ is an $r\times r$ very badly approximable unitary-valued function such that $\|H_U\|_{\rm e}<1$, and
$\Psi$ is an $(m-r)\times(n-r)$ matrix function such that  $\|\Psi\|_{L^\be}\le\s$, $\|H_\Psi\|<\s$,
and $\|H_\Psi\|_{\rm e}\le\|H_\Phi\|_{\rm e}$.

Moreover, $\Phi$ is very badly approximable if and only if $\Psi$ is very badly approximable.}

\medskip

{\bf Remark 1.}
If $m=r$ or $n=r$, by $\left(\begin{array}{cc}\s U&\bs{0}\\\bs{0}&\Psi\end{array}\right)$
we mean $\left(\begin{array}{cc}\s U&\bs{0}\end{array}\right)$ or 
$\left(\begin{array}{cc}\s U\\\bs{0}\end{array}\right)$ respectively, in which case $\Phi$ is very badly
approximable if and only if $\Phi$ is badly approximable.

\medskip

Such factorizations are a special case of partial canonical factorizations. 
Partial canonical factorizations in the general case are defined in \cite{AP}.

\medskip

{\bf Remark 2.}
Actually, if $\Phi$ admits a factorization as above, then $\Phi$ must be badly approximable even without
the assumption  $\|H_\Phi\|_{\rm e}<\|H_\Phi\|$.

\medskip

{\bf Remark 3.} Note that if $U$ is a very badly approximable unitary-valued function such that
$\|H_U\|_{\rm e}<1$, then the Toeplitz operator $T_U$ is Fredholm, see \cite{AP}.

\medskip

Let us now define  a canonical factorization. Let $\s_0,\cdots,\s_{\iota-1}$ 
be all distinct nonzero superoptimal singular values of $\Phi$. Suppose that $d_j$ is the multiplicity 
of the superoptimal singular value $\s_j$ of $\Phi$. A {\it canonical factorization} of $\Phi$ is a representation of
$\Phi$ of the form
\bay
\label{cf}
\Phi=\W_0^*\cdots\W^*_{\iota-1}
\left(\begin{array}{ccccc}\s_0U_0&\0&\cdots&\0&\bs{0}\\
\0&\s_1U_1&\cdots&\0&\0\\
\vdots&\vdots&\ddots&\vdots&\vdots\\
\0&\0&\cdots&\s_{\iota-1}U_{\iota-1}&\0\\
\0&\0&\cdots&\0&\0\end{array}\right)
\V^*_{\iota-1}\cdots\V_0^*,
\ey
where the $U_j$ are $d_j\times d_j$
unitary-valued very badly approximable matrix functions such that $\|H_{U_j}\|_{\rm e}<1$,
the matrix functions $\V_j$ and $\W_j$, $1\le j\le\iota-1$, have the form
$$
\V_j=\left(\begin{array}{cc}\bs{I}_{d_0+\cdots+d_{r-1}}&\0\\\0&\breve{\V}_j\end{array}\right)
\quad\mbox{and}\quad
\W_j=\left(\begin{array}{cc}\bs{I}_{d_0+\cdots+d_{r-1}}&\0\\\0&\breve{\W}_j\end{array}\right),
\quad1\le j\le\iota-1,
$$
$\V_0$ and $\W_0^{\rm t}$ are $d_0$-balanced matrix functions and $\breve{\V}_j$ and
$\breve{\W}_j^{\rm t}$ are $d_j$-balanced matrix functions.
Note that the last zero row has size $(m-(d_0+\cdots+d_{\iota-1}))\times n$.
If $m=d_0+\cdots+d_{\iota-1}$, this means that there is no zero row in \rf{cf}. 
A similar remark can be made about the last zero column in \rf{cf}.

It was shown in \cite{AP} that an admissible matrix function
$\Phi$ is very badly approximable if and only if it admits a canonical factorization.
Again, if $\Phi$ is an arbitrary bounded matrix function (not necessarily admissible) that
admits a canonical factorization, then $\Phi$ must be very badly approximable.

Finally, we need the following result from \cite{AP}.

\medskip

{\bf Theorem E.}  {\it Let  $U$ be a unitary-valued matrix function such that 
$\|H_U\|_{\rm e}<1$. Then $U$ is very badly approximable if and only if the Toeplitz operator
$T_{\bar zU^*}$ has trivial kernel.}

\medskip

Note that all the above results can be found in Chapter 14 of the book \cite{P}.

\

\section{\bf Analytic Families of Subspaces}
\label{s3}
\setcounter{equation}{0}

\

In this section we are going to state one more necessary condition for an admissible matrix function to be
very badly approximable. This condition involves analyticity of certain families of subspaces.

Let $\Phi$ be a matrix function in $L^\be(\mm_{m,n})$ and let $\s>0$.
%such that the functions
%$\z\mapsto s_j(\z)$ are constant almost everywhere on $\T$ and suppose that for almost all $\z\in\T$
%the numbers $\s_0,\cdots,\s_{\iota-1}$ are all distinct nonzero singular values of $\Phi(\z)$, which are 
%arranged in the decreasing order.
For $\z\in\T$ we denote by ${\frak S}_\Phi^{(\s)}(\z)$
the linear span of all Schmidt vectors \footnote{Recall that if $A$ is 
an $m\times n$ matrix and $s$ is a singular value of $A$, 
a nonzero vector $x\in\C^n$ is called
a {\it Schmidt vector} corresponding to $s$ if $A^*Ax=s^2x$.}
 of $\Phi(\z)$ that correspond to the singular values of $\Phi(\z)$ that are greater than or equal to $\s$.
The subspaces ${\frak S}_\Phi^{(\s)}(\z)$ are defined for almost all $\z\in\T$. 

As we have mentioned in the introduction, in \cite{PY2} an example of a continuous $2\times2$ matrix function
$\Phi$ was given such that $T_{z\Phi}$ is invertible, $s_0(\Phi(\z))=1$, $s_1(\Phi(\z))=\a$, $\a\in(0,1)$,
but $\Phi$ is not badly approximable. If we look at the subspace of maximizing vectors of $\Phi(\z)$, $\z\in\T$,
in that example, we can easily observe that the family of subspaces ${\frak S}_\Phi^{(1)}(\z)$, $\z\in\T$, 
is not analytic in the following sense.

\medskip

{\bf Definition.} Let ${\frak L}_n$ be the set of all subspaces of $\C^n$.
A family of subspaces $L(\z)$, $\z\in\T$, defined for almost all $\z\in\T$ is called {\it analytic} if there exist 
functions $\xi_1,\cdots,\xi_k$ in $H^2(\C^n)$ such that $L(\z)=\spn\{\xi_j(\z):~1\le j\le k\}$ for almost 
all $\z\in\T$.

\medskip

{\bf Remark 1.} It is easy to see that if
$L(\z)$, $\z\in\T$, is an analytic family of subspaces, then there exists $r\in\Z_+$ such that
$\dim L(\z)=r$ everywhere on $\T$ and there exist $\xi_1,\cdots,\xi_r$ in $H^2(\C^n)$ such that
$L(\z)=\spn\{\xi_j(\z):~1\le j\le r\}$ for almost all $\z\in\T$.

\medskip

In the Introduction we have mentioned the following necessary conditions for an admissible matrix
function $\Phi$ to be very badly approximable:

\medskip

\begin{enumerate}
\item[(C1)] {\it the functions $\z\mapsto s_j(\Phi(\z))$, $0\le j\le\min\{m,n\}-1$,
 are constant almost everywhere on $\T$};

\item[(C2)] $\Ker T_{\bar z\Phi^*}=\{f\!\in\!H^2(\C^n):\Phi^*f=\0\}$
{\it and} $\Ker T_{\bar z\ov{\Phi}}=\{f\!\in\! H^2(\C^n):\ov{\Phi}f=\0\}$.
\end{enumerate}

\medskip

In this section we consider the following important condition:

\begin{enumerate}
\item[{\rm(C3)}] {\it if $\s>0$, then ${\frak S}_\Phi^{(\s)}(\z)$, $\z\in\T$, 
and ${\frak S}_{\Phi^{\rm t}}^{(\s)}(\z)$, $\z\in\T$, are analytic families of subspaces.}
\end{enumerate}

\begin{thm}
\label{1mo}
Let $\Phi$ be an admissible very badly approximable matrix function in $L^\be(\mm_{m,n})$.
Then $\Phi$ satisfies {\em(C3)}.
\end{thm}

We will see later that Theorem \ref{1mo} is an immediate consequence of Theorem \ref{plo}.

\medskip

{\bf Remark 2.} Note that it follows easily from the above Remark 1 that the analyticity of
the families
${\frak S}_\Phi^{(\s)}(\z)$, $\z\in\T$, for $\s>0$ implies condition (C1). A fortiori (C3) 
implies (C1). 

Indeed, for any $\s>0$ the analytic family of subspaces ${\frak S}_\Phi^{(\s)}$ has constant
dimension a.e.~on $\T$, and as one can easily see, this is possible only if  the functions
$\z\mapsto s_j(\Phi(\z))$ are constant almost everywhere on $\T$. 

\medskip

We show in this section that if $\Phi$ is an admissible matrix function satisfying (C3),
then $\Phi$ admits a factorization of the form \rf{cf} with $\V_j$ and $\W_j$ as in \rf{cf} and
unitary-valued functions $U_j$ such that $\|H_{U_j}\|_{\rm e}<1$. We call such factorizations {\it quasicanonical}.
(A quasicanonical factorization is canonical if the unitary-valued functions $U_j$ are very badly approximable).

Then we show that conditions 
(C1)--(C3) are not sufficient for an admissible function $\Phi$ to be very badly approximable.

Note here that the condition that the families ${\frak S}_\Phi^{(\s)}(\z)$, $\z\in\T$,  are analytic for $\s>0$
does not imply that
the families ${\frak S}_{\Phi^{\rm t}}^{(\s)}(\z)$, $\z\in\T$ are analytic for $\s>0$ (even under condition (C2))
as the following example shows.

\medskip

{\bf Example 1.} Let 
$$
W=\left(\begin{array}{cc}w_1&-\ov{w}_2\\w_2&\ov{w}_1\end{array}\right)
$$
be a thematic (1-balanced) matrix function, i.e., 
$w_1,\,w_2\in H^\be$, $|w_1|^2+|w_2|^2=\bs{1}$, and $w_1$ and $w_2$ are coprime.
Consider the function 
$$
\Phi=W^*\left(\begin{array}{cc}\bar z&\bs{0}\\\bs{0}&\frac{\bar z}2\end{array}\right)
$$
Clearly, ${\frak S}_\Phi^{(\s)}$ is a constant function for each $\s>0$, and so the family 
${\frak S}_\Phi^{(\s)}$, $\z\in\T$,
is analytic. 

Let us verify that $\Phi$ satisfies (C2). Suppose that 
$g\in\Ker T_{\bar z\Phi^*}$. Clearly, this means that
$Wg\in H^2_-(\C^2)$, i.e., $g\in\Ker T_W$. By Theorem A in \S \ref{s2}, $g=\0$. 
Similarly, it is easy to see that $\Ker T_{\bar z\ov{\Phi}}=\{\0\}$ if and only is
$\Ker T_{W^{\rm t}}=\{\0\}$. The last equality also follows from Theorem A.

Let us show that the family ${\frak S}_{\Phi^{\rm t}}^{(\s)}(\z)$, $\z\in\T$, does not have to be analytic. Suppose that 
$f=\left(\begin{array}{c}f_1\\f_2\end{array}\right)\in H^2(\C^2)$ and $f(\z)$
is a maximizing vector of $\Phi^{\rm t}(\z)$ for almost all $\z\in\dd$. Clearly,
$\ov{W}f$ must be of the form 
$$
\ov{W}f=\left(\begin{array}{c}\f\\\0\end{array}\right),\quad \f\in L^2.
$$
Since $W$ is a unitary-valued matrix function, it follows that 
$$
f=W^{\rm t}\ov{W}f=W^{\rm t}\left(\begin{array}{c}\f\\\0\end{array}\right)=
\left(\begin{array}{c}\f w_1\\-\f\ov{w}_2\end{array}\right).
$$
Thus the function ${\frak S}_{\Phi^{\rm t}}^{(1)}$ is analytic if and only if there exists a function
$\f\in L^2$ such that $\f\w_1\in H^2$ and $\f\ov{w}_2\in H^2$. Suppose now that $w_1$ is invertible in $H^\be$. Then
$\f$ must be in $H^2$. Then the function $w_2$ must have a meromorphic pseudocontinuation (see \cite{N}, Lect. II, Sect. 1).
Hence, if $w_1$ is invertible in $H^\be$ and $w_2$ does not have a pseudocontinuation, the function 
${\frak S}_{\Phi^{\rm t}}^{(1)}$ is not analytic. $\bl$

\medskip

The following example shows that none of the two conditions in (C2) implies the other one (even under conditions 
(C1) and (C3)).

\medskip

{\bf Example 2.} Let $V=\left(\begin{array}{cc}v_1&-\bar{v}_2\\v_2&\bar{v}_1\end{array}\right)$ be a continuous 
thematic ($1$-balanced) matrix function. Consider the matrix function $\Phi$ defined by
$$
\Phi=\left(\begin{array}{cc}\1&\0\\\0&\0\end{array}\right)V^*=
\left(\begin{array}{cc}\bar v_1&\bar v_2\\\0&\0\end{array}\right).
$$
Obviously, $\Phi$ satisfies (C1) and (C3). Let us show that $\Phi$ satisfies the first condition in (C2). Suppose that 
$\left(\begin{array}{c}g_1\\g_2\end{array}\right)\in\Ker T_{\bar z\Phi^*}$. Then
$$
\bar z\Phi^*\left(\begin{array}{c}g_1\\g_2\end{array}\right)=
\left(\begin{array}{c}\bar zv_1g_1\\\bar zv_2g_1\end{array}\right)\in H^2_-(\C^2).
$$
It follows that both $v_1g_1$ and $v_2g_1$ are constant functions. Suppose now that both $v_1$ and $v_2$ 
are nonzero functions such that the function $v_1v_2^{-1}$ is nonconstant. It is easy to see that in this case
$g_1=\0$. Thus
$$
\Phi^*\left(\begin{array}{c}g_1\\g_2\end{array}\right)=\Phi^*\left(\begin{array}{c}\0\\g_2\end{array}\right)=
\left(\begin{array}{c}\0\\\0\end{array}\right).
$$
However, $\Ker T_{\bar z\ov{\Phi}}\ne\{f\!\in\! H^2(\C^n):\ov{\Phi}f=\0\}$. Indeed, 
$$
\bar z\ov{\Phi}\left(\begin{array}{c}g_1\\g_2\end{array}\right)=
\bar z\left(\begin{array}{c}v_1g_1+v_2g_2\\\0\end{array}\right).
$$
Clearly, we can choose nonzero functions $g_1$ and $g_2$ in $H^2$ such that $v_1g_1+v_2g_2=\1$. $\bl$

\begin{thm}
\label{scf}
Let $\Phi$ be a matrix function in $L^\be(\mm_{m,n})$ that satisfies {\em(C3)}.
Then $\Phi$ admits a quasicanonical factorization.
\end{thm}

\Pf As we have already observed, (C3) implies (C1).  Let $\s_0>\cdots>\s_{\iota-1}$ 
be positive numbers such that for almost all $\z\in\T$ the distinct nonzero singular values of $\Phi(\z)$ 
are precisely, $\s_0,\cdots,\s_{\iota-1}$. We argue by induction on $\iota$. If $\iota=0$, then $\Phi=\bs{0}$.

Let now $\iota>0$. Suppose that $\dim{\frak S}_\Phi^{(\s_0)}(\z)=r$ for almost all $\z\in\T$. Obviously,
$\dim{\frak S}_{\Phi^{\rm t}}^{(\s_0)}(\z)=r$ for almost all $\z\in\T$. 
Let us first show that $\Phi$ admits a factorization of the form
\bay
\label{spc}
\Phi=\W^*\left(\begin{array}{cc}\s_0 U&\bs{0}\\\bs{0}&\Psi\end{array}\right)\V^*,
\ey 
in which $\V$ and $\W^{\rm t}$ 
are $r$-balanced unitary-valued matrix functions, $U$ is an $r\times r$ 
unitary-valued matrix function such that $\|H_U\|_{\rm e}<1$. The proof is similar to
the proof of Theorem 4.3 of \cite{AP}.

Let $\xi_1,\cdots,\xi_r$ and
$\eta_1,\cdots,\eta_r$ are functions in $H^2(\C^r)$ such that
$$
{\frak S}_\Phi^{(\s_0)}(\z)=\spn\{\xi_1(\z),\cdots,\xi_r(\z)\}\quad\mbox{and}\quad
{\frak S}_{\Phi^{\rm t}}^{(\s_0)}(\z)=\spn\{\eta_1(\z),\cdots,\eta_r(\z)\}
$$
almost everywhere on $\T$. Let $\cL$ be the minimal invariant subspace of multiplication by $z$ on $H^2(\C^n)$
that contains $\xi_1,\cdots,\xi_r$ and let $\M$ be the minimal invariant subspace of multiplication 
by $z$ on $H^2(\C^n)$ that contains $\eta_1,\cdots,\eta_r$. 

It is easy to see from \rf{r} that there exist $n\times r$
inner functions $\U$ and $\O$ such that $\cL=\U H^2(\C^r)$ and $\M=\O H^2(\C^r)$.
Let us show that $\U$ and $\O$ are co-outer.

Indeed, suppose that $\U^{\rm t}=\L F$, where $\L$ is an inner matrix function and $F$ is an outer matrix function. 
Since $\dim {\frak S}_\Phi^{(\s_0)}(\z)=r$ for almost all $\z\in\T$, it follows that $\rank \L(\z)=r$ almost everywhere
on $\T$, and so $\L$ is an $r\times r$ inner function, and so $F^{\rm t}$ is inner.
Since $\U=F^{\rm t}\L^{\rm t}$, it follows that $\cL=\U H^2(\C^r)\subset F^{\rm t}H^2(\C^r)$.
Clearly, for every $d\in\C^r$,
the vector $\U(\z)d$ belongs to ${\frak S}_\Phi^{(\s_0)}(\z)$ for almost all $\z\in\T$. It follows that 
$F^{\rm t}(\z)d=\U(\z)\ov{\L(\z)}d\in{\frak S}_\Phi^{(\s_0)}(\z)$ for almost all $\z\in\T$, and so 
$F^{\rm t}H^2(\C^r)=\cL=\U H^2(\C^r)$. Hence, $\L$ is a constant matrix (see \cite{N}) and $\U$ is co-outer.

Let now  $\Theta$ and $\Xi$ be inner and co-outer matrix 
functions such that the matrix functions
$$
\V=\left(\begin{array}{cc}\U&\ov{\Theta}\end{array}\right)\quad\mbox{and}\quad
\W^{\rm t}=\left(\begin{array}{cc}\O&\ov{\Xi}\end{array}\right)
$$
are $r$-balanced.

It is easy to see that if $q_1,\cdots,q_r$ are scalar polynomials and $\xi=q_1\xi_1+\cdots+q_r\xi_r$,
then $\xi(\z)$ is a maximizing vector of $\Phi(\z)$ for almost all $\z\in\T$. It follows that for any
function $f\in\cL$ the vector $f(\z)$ is a maximizing vector of $\Phi(\z)$ for almost every $\z\in\T$.
In particular, the columns of $\U(\z)$ are maximizing vectors of $\Phi(\z)$ almost everywhere on $\T$.
For the same reason, the columns of $\O(\z)$ are maximizing vectors of $\Phi^{\rm t}(\z)$ for almost 
every $\z\in\T$.

We need two obvious and well known lemmas.

\begin{lem}
\label{ol1}
Let $A\in\mm_{m,n}$ and $\|A\|=1$. Suppose that $v_1,\cdots,v_r$ is an 
orthonormal family of maximizing vectors of $A$ and $w_1,\cdots,w_r$ is an 
orthonormal family of maximizing vectors of $A^{\text t}$. Then{\em
$$
\left(\begin{array}{ccc}w_1&\cdots&w_r\end{array}\right)^{\text t}A
\left(\begin{array}{ccc}v_1&\cdots&v_r\end{array}\right)
$$}
is a unitary matrix.
\end{lem}

\begin{lem}
\label{ol2}
Let $A$ be a matrix in $\mm_{m,n}$ such that $\|A\|=1$ and $B$ has the form
$$
A=\left(\begin{array}{cc}A_{11}&A_{12}\\A_{21}&A_{22}\end{array}\right),
$$
where $A_{11}$ is a unitary matrix. Then $A_{12}$ and $A_{21}$ are the zero matrices.
\end{lem}

Consider the matrix function
$$
\left(\begin{array}{cc}U&X\\Y&\Psi\end{array}\right)\df\s_0^{-1}\W\Phi\V.
$$
It follows easily from Lemmas \ref{ol1} and \ref{ol2} that 
$U=\s_0^{-1}\O^{\rm t}\Phi\U$ is a unitary-valued matrix function while
$X$ and $Y$ are the zero matrix functions. Thus
\rf{spc} holds.

Let us show that $\|H_U\|_{\rm e}<1$. We have
\bey
\|H_U\|_{\text e}&=&\dist_{L^\be}\big(U,(H^\be+C)(\mm_{r,r})\big)\\
&=&\s_0^{-1}\dist_{L^\be}\big(\O^{\text t}\Phi\U,(H^\be+C)(\mm_{r,r})\big)\\
&\le&\s_0^{-1}\dist_{L^\be}\big(\Phi,(H^\be+C)(\mm_{m,n})\big)=\s_0^{-1}\|H_\Phi\|_{\text e}<1.
\eey

It is sufficient to show that $\Psi$ has a quasicanonical factorization.
Clearly, for almost all $\z\in\T$, $\s_1,\cdots,\s_{\iota-1}$ are all distinct nonzero singular values of 
$\Psi(\z)$. If $\iota=1$, then $\Psi=\bs{0}$, and everything is trivial.
Let us show that the families ${\frak S}_\Psi^{(\s_j)}(\z)$, $\z\in\T$, are analytic for $j\ge1$.

Consider the family ${\frak S}_\Phi^{(\s_{j})}(\z)$, $\z\in\T$. Let $\xi_1,\cdots,\xi_\kappa$ be functions in $H^2(\C^n)$ 
such that 
$$
{\frak S}_\Phi^{(\s_{j})}(\z)=\spn\{\xi_1(\z),\cdots,\xi_\kappa(\z)\}\quad\mbox{for almost all}\quad\z\in\T.
$$
Since $\V$ is unitary-valued, we have
$$
{\frak S}_{\Phi\V}^{(\s_j)}(\z)=\spn\{(\V^*\xi_1)(\z),\cdots,(\V^*\xi_\kappa)(\z)\}
\quad\mbox{for almost all}\quad\z\in\T.
$$
We have $\V^*\xi_j=\left(\begin{array}{c}\U^*\xi_r\\\Theta^{\rm t}\xi_r\end{array}\right)$. Obviously,
$\left(\begin{array}{c}(\U^*\xi_r)(\z)\\0\end{array}\right)$ is a maximizing vector of
$(\Phi\V)(\z)$, and so it belongs to ${\frak S}_{\Phi\V}^{(\s_j)}(\z)$. Thus 
$\left(\begin{array}{c}0\\(\Theta^{\rm t}\xi_r)(\z)\end{array}\right)\in{\frak S}_{\Phi\V}^{(\s_j)}(\z)$. 
It is easy to see that $(\Theta^{\rm t}\xi_r)(\z)\in{\frak S}_{\Psi}^{(\s_j)}(\z)$. Moreover, it is evident that
$$
{\frak S}_{\Psi}^{(\s_j)}(\z)=\spn\{(\Theta^{\rm t}\xi_r)(\z):~1\le r\le\kappa\}\quad
\mbox{for almost all}\quad\z\in\T,
$$
which proves that ${\frak S}_{\Psi}^{(\s_j)}(\z)$, $\z\in\T$, is an analytic
family of subspaces. The same reasoning shows that the functions
${\frak S}_{\Psi^{\rm t}}^{(\s)}(\z)$ $\z\in\T$, is an analytic
family of subspaces. for $\s>0$.

By the inductive hypothesis, $\Psi$ admits a quasicanonical factorization. $\bl$

It turns out, however, that conditions (C1)--(C3) do not imply that the matrix function $\Phi$ is very badly 
approximable.

\medskip

{\bf Example 3.} Consider the function
$$
\Phi=\left(\begin{array}{cc}\1&\0\\\0&\frac12\bar z^2\end{array}\right)
\left(\begin{array}{cc}\frac1{\sqrt{2}}\bar z&\frac1{\sqrt{2}}\\[.2cm]-\frac1{\sqrt{2}}&\frac1{\sqrt{2}}z\end{array}\right),
$$
Let us show that $\Phi$ satisfies (C1)--(C3), but $\Phi$ is not even badly approximable. Note that
$$
\bar z\Phi=\left(\begin{array}{cc}\bar z&\0\\\0&\frac12\bar z^3\end{array}\right)
\left(\begin{array}{cc}\frac1{\sqrt{2}}\bar z&\frac1{\sqrt{2}}\\[.2cm]-\frac1{\sqrt{2}}&\frac1{\sqrt{2}}z\end{array}\right)
$$
is a canonical factorization of $\bar z\Phi$, and so $\bar z\Phi$ is very badly approximable 
(see \S 
\ref{s2}). Hence, it satisfies (C1)--(C3). Clearly, conditions (C1) and (C3) are invariant
under multiplication by $z$. Thus
$\Phi$ satisfies (C1) and (C3).

Let us show that $\Ker T_{\bar z\Phi^*}=\{\0\}$. Suppose that 
$\left(\begin{array}{c}g_1\\g_2\end{array}\right)\in\Ker T_{\bar z\Phi^*}$. We have 
$$
\bar z\Phi^*\left(\begin{array}{c}g_1\\g_2\end{array}\right)=
\frac1{2\sqrt{2}}\left(\begin{array}{c}2g_1-zg_2\\2\bar zg_1+g_2\end{array}\right).
$$
Thus $2g_1-zg_2\in H^2_-$ and $2\bar zg_1+g_2\in H^2_-$. Multiplying the first inclusion by $\bar z$, we obtain
$2\bar zg_1-g_2\in H^2_-$, and so
$\bar zg_1\in H^2_-$ and $g_2\in H^2_-$. This implies that $g_2=\0$, and it follows from the first inclusion that $g_1=\0$.

Let us prove now that $\Ker T_{\bar z\ov{\Phi}}=\{\0\}$. Suppose that 
$\left(\begin{array}{c}g_1\\g_2\end{array}\right)\in\Ker T_{\bar z\ov{\Phi}}$. We have 
$$
\bar z\ov{\Phi}\left(\begin{array}{c}g_1\\g_2\end{array}\right)=
\frac1{2\sqrt{2}}\left(\begin{array}{c}2g_1+2\bar zg_2\\-zg_1+g_2\end{array}\right).
$$
It follows that $g_1+\bar zg_2\in H^2_-$ and $-zg_1+g_2\in H^2_-$. Again, 
multiplying the second inclusion by $\bar z$, we obtain $g_1-\bar zg_2\in H^2_-$, and so both $g_1$ and $\bar z g_2$ belong
to $H^2_-$. Thus $g_1=\0$, and it follows from the second inclusion that $g_2=\0$.

We can show now that $\Phi$ is not even badly approximable. Clearly, $\|\Phi\|_{L^\be}=1$. 
If $\Phi$ is badly approximable, then
$\|H_\Phi\|=1$. Since $\Phi$ is continuous, $H_\Phi$ is compact and so $H_\Phi$ has a maximizing vector 
$f=\left(\begin{array}{c}f_1\\f_2\end{array}\right)$. 
Put
$$
V=
\left(\begin{array}{cc}\frac1{\sqrt{2}}\bar z&\frac1{\sqrt{2}}\\[.2cm]
-\frac1{\sqrt{2}}&\frac1{\sqrt{2}}z\end{array}\right)^*
\df\left(\begin{array}{cc}v_1&-\bar v_2\\v_2&\bar v_1\end{array}\right).
$$
Clearly, 
the second component of the vector function $V^*f$
must be zero
and $\Phi f$ must belong to $H^2_-(\C^2)$. Let $V^*f=\left(\begin{array}{c}h\\\0\end{array}\right)$,
where $h\in L^2$.
We have
$$
f=VV^*f=V\left(\begin{array}{c}h\\\0\end{array}\right)=\left(\begin{array}{c}v_1h\\v_2h\end{array}\right)\in H^2(\C^2).
$$
Since the matrix function $\left(\begin{array}{c}v_1\\v_2\end{array}\right)$ is co-outer, it is easy to see that $h\in H^2$.
We have 
$$
\Phi f=\left(\begin{array}{cc}\1&\0\\\0&\frac12\bar z^2\end{array}\right)V^*f
=\left(\begin{array}{cc}\1&\0\\\0&\frac12\bar z^2\end{array}\right)\left(\begin{array}{c}h\\\0\end{array}\right)
=\left(\begin{array}{c}h\\\0\end{array}\right)\in H^2_-(\C^2),
$$
and so $h=\0$. Hence, $H_\Phi$ has no maximizing vector and we get a contradiction. $\bl$

\

\section{\bf Very Badly Approximable Matrix Functions}
\label{s4}
\setcounter{equation}{0}

\

We obtain in this section a necessary and sufficient condition for an admissible matrix function
to be very badly approximable. To do this, we slightly modify the necessary conditions stated in the previous
section.

\medskip

{\bf Definition.} Let ${\frak L}_n$ be the set of all subspaces of $\C^n$.
Suppose that $L:\T\to{\frak L}_n$ is an ${\frak L}_n$-valued function defined almost everywhere.
We say that functions $\xi_1,\cdots,\xi_l$ in $H^2(\C^n)$ {\it span the function} $L$ if
$L(\z)=\spn\{\xi_j(\z):~1\le j\le l\}$ for almost all $\z\in\T$.

\medskip

It is easy to see that if functions $\xi_1,\cdots,\xi_l$ in $H^2(\C^n)$ span an ${\frak L}_n$-valued function $L$,
then $\dim L(\z)$ is constant for almost all \mbox{$\z\in\T$} and there  exist \lb functions
$\eta_1,\cdots,\eta_k$ in \mbox{$\spn\{\xi_j:~1\le j\le l\}$} such that \mbox{$k=\dim L(\z)$} and
\lb$L(\z)=\spn\{\eta_j(\z):~1\le j\le k\}$ almost everywhere on $\T$.

As in \S \ref{s3}, we consider a matrix function
$\Phi$ in $L^\be(\mm_{m,n})$ and for $\s>0$ we associate with $\Phi$ the linear span 
${\frak S}_\Phi^{(\s)}(\z)$ 
 of all Schmidt vectors of $\Phi(\z)$ that correspond to the singular values greater than or equal to $\s$.

We consider in this section the following condition:

\begin{enumerate}
\item[(C4)] {\em for each $\s>0$, the analytic family of subspaces ${\frak S}_\Phi^{(\s)}$
is spanned by finitely many functions in
$\Ker T_\Phi$}.
\end{enumerate}

Recall that a function $\Phi\in L^\be(\mm_{m,n})$ is called {\em admissible} if 
$t_k(\Phi)>\|\Phi\|_{\rm e}$ for all nonzero superoptimal values  $t_k(\Phi)$. In
particular, any continuous matrix-valued function is admissible, since the Hankel operator
$H_\Phi$ is compact (and so its essential norm is $0$) in this case. 

\begin{thm}
\label{plo}
If $\Phi$ is an admissible very badly approximable matrix function in $L^\be(\mm_{m,n})$,
then $\Phi$ satisfies {\em(C4)}.

Conversely, if $\Phi$ is an arbitrary function in 
$L^\be(\mm_{m,n})$ that satisfies  {\em(C4)}, then $\Phi$ is very badly approximable
and $\0$ is the only superoptimal approximant of $\Phi$. 
\end{thm}

{\bf Remark 1.} Clearly, condition (C4) implies that ${\frak S}_\Phi^{(\s)}(\z)$, $\z\in\T$, 
is an analytic family of subspaces, and it is easy to see that 
Theorem \ref{plo} implies Theorem \ref{1mo}.

\medskip

{\bf Remark 2.} As we have already observed (see Remark 2 after Theorem \ref{1mo}), 
condition (C4) implies that the functions $\z \mapsto s_j(\Phi(\z))$ are constant
almost everywhere on $\T$. 
\medskip

{\bf Remark 3.} It is interesting to observe that to prove that (C4) implies that $\Phi$ is
very badly approximable, we do not need the fact that $\Phi^{\rm t}$ satisfies (C4).

\medskip

{\bf Proof of Theorem \ref{plo}.}  
Suppose first that $\Phi$ is admissible and very badly approximable.
Then $s_j(\Phi(\z))=t_j(\Phi)$, $0\le j\le\min\{m,n\}-1$, almost everywhere on $\T$ (see \rf{sv}).
Let us prove by induction on
$\vk$ that if $\Phi$ is an admissible very badly approximable matrix function and 
for almost all $\z\in\T$, 
$$
\s_0=\s_0(\Phi)>\s_1=\s_1(\Phi)>\s_2=\s_2(\Phi)>\cdots
$$
are all distinct nonzero singular values of $\Phi(\z)$, then
${\frak S}_\Phi^{(\s_\vk)}={\frak S}_\Phi^{(\s_\vk(\Phi))}$ is spanned by finitely many functions in $\Ker
T_\Phi$.  

By Theorem D stated in \S 2,
$\Phi$ admits a factorization
\bay
\label{pcf}
\Phi=\W^*\left(\begin{array}{cc}\s_0 U&\bs{0}\\\bs{0}&\Psi\end{array}\right)\V^*,
\ey
where $\s_0=\|\Phi\|_{L^\be}$,
$\V$ and $\W^{\rm t}$ are $r$-balanced unitary-valued functions, 
\lb$1\le r\le\min\{m,n\}$, $U$ is an $r\times r$ very badly approximable unitary-valued function such that 
$T_U$ is Fredholm and $\|H_U\|_{\rm e}<1$, and $\Psi$ is an admissible very badly approximable
matrix function with $\|\Psi\|_{L^\be}=\s_1<\s_0$. 

Let us prove  first that ${\frak S}_\Phi^{(\s_0)}$ is spanned by finitely many functions in $\Ker T_\Phi$.
Since $T_U$ is Fredholm (see Remark 3 after Theorem D in \S 2),  it admits a Wiener--Hopf
factorization
$$
U=G^*DF,
$$
where
$$
D=\left(\begin{array}{cccc}z^{d_1}&\bs{0}&\cdots&\bs{0}\\
\bs{0}&z^{d_2}&\cdots&\bs{0}\\
\vdots&\vdots&\ddots&\vdots\\
\bs{0}&\bs{0}&\cdots&z^{d_r}
\end{array}\right),
$$ 
$F$ and $G$ are $r\times r$ matrix functions such that $F^{\pm1}\in H^2(\mm_{r,r})$ and 
\lb$G^{\pm1}\in H^2(\mm_{r,r})$, and $d_1,\cdots,d_r\in\Z$ 
(Simonenko's theorem; see e.g., \cite{P}, Ch. 3, \S 5). 
By Theorem E, $\Ker T_{\bar zU^*}=\{\0\}$,
which implies easily that the indices $d_1,\cdots,d_r$ are negative. Let $c_1,\cdots,c_r$ be a basis in $\C^r$.
Consider the functions
\bay
\label{spn}
\U^*F^{-1}\bs{c}_1,\cdots,\U^*F^{-1}\bs{c}_r,
\ey
where $\bs{c}_j$ denotes the constant function identically equal to $c_j$. Since $\V$ is a unitary-valued
function, it is easy to see that 
$$
\V^*\U^*F^{-1}\bs{c}_j=\left(\begin{array}{c}F^{-1}\bs{c}_j\\\bs{0}\end{array}\right),
$$
and so 
$$
\Phi\U^*F^{-1}\bs{c}_j=\W^*\left(\begin{array}{cc} G^*D\bs{c}_j\\\bs{0}
\end{array}\right)=
\ov{\O}G^*D\bs{c}_j\in H^2_-(\C^m),
$$
since the Wiener--Hopf indices $d_j$ are negative. It is easy to see now that the functions in \rf{spn}
belong to $\Ker T_\Phi$ and span ${\frak S}_\Phi^{(\s_0)}$.

Let now that $\vk>0$. Clearly, for almost all $\z\in\T$,
$\s_1>\s_2\cdots>$  are all nonzero singular values of $\Psi(\z)$ and $\s_\vk=\s_\vk(\Phi)=\s_{\vk-1}(\Psi)$.
By the inductive hypothesis, there exist functions $\xi_1,\cdots,\xi_l$
in $\Ker H_\Psi$ that span ${\frak S}_\Psi^{(\s_{\vk})}={\frak S}_\Psi^{(\s_{\vk-1}(\Psi))}$. 
By Theorem C, the functions $\Theta$ and $\Xi$ are left invertible
in $H^\be$. Let $Q\in H^\be(\mm_{n-r,n})$ and $R\in H^\be(\mm_{m-r,m})$ such that
$Q\Theta=\bs{I}_{n-r}$ and $R\Xi=\bs{I}_{m-r}$. Put
$$
\eta_j=Q^{\rm t}\xi_j+\U q_j,\quad 1\le j\le l.
$$
where the functions $q_j\in H^2(\C^r)$ will be chosen later. We have
\bay
\label{Veta}
\V^*\eta_j=\left(\begin{array}{c}\U^*\\\Theta^{\rm t}\end{array}\right)(Q^{\rm t}\xi_j+\U q_j)=
\left(\begin{array}{c}\U^*Q^{\rm t}\xi_j+q_j\\\xi_j\end{array}\right),
\ey
since $\V$ is unitary-valued and $\Theta^{\rm t}Q^{\rm t}=\bs{I}_{n-r}$.

Since $\W$ is a unitary-valued function, we obtain
$$
\bs{I}_m=\W^*\W=\ov{\O}\O^{\rm t}+\Xi\Xi^*,
$$
and so
\bay
\label{Xi}
\Xi=\Xi(R\Xi)^*=\Xi\Xi^*R^*=(\bs{I}_m-\ov{\O}\O^{\rm t})R^*.
\ey 
We have now from \rf{Veta} and \rf{Xi}
\bey
\Phi\eta_j&=&\W^*\left(\begin{array}{cc}\s U&\bs{0}\\\bs{0}&\Psi\end{array}\right)
\left(\begin{array}{c}\U^*Q^{\rm t}\xi_j+q_j\\\xi_j\end{array}\right)\\
&=&\W^*\left(\begin{array}{c}\s U(\U^*Q^{\rm t}\xi_j+q_j)\\\Psi\xi_j\end{array}\right)\\
&=&\s\ov{\O}U(\U^*Q^{\rm t}\xi_j+q_j)+(\bs{I}_m-\ov{\O}\O^{\rm t})R^*\Psi\xi_j\\
&=&R^*\Psi\xi_j+\ov{\O}\big(\s U(\U^*Q^{\rm t}\xi_j+q_j)-\O^{\rm t}R^*\Psi\xi_j\big).
\eey
In order that $\Phi\eta_j\in H^2_-(\C^m)$, it is sufficient that 
$$
\pp_+\big(\s U(\U^*Q^{\rm t}\xi_j+q_j)-\O^{\rm t}R^*\Psi\xi_j\big)=\bs{0},
$$
which means that
$$
\s T_Uq_j=\pp_+\big(\O^{\rm t}R^*\Psi\xi_j-\s U\U^*Q^{\rm t}\xi_j\big).
$$
Since $\Range T_U=H^2(\C^r)$, we can find a solution $q_j\in H^2(\C^r)$. This proves that
$\eta_j\in \Ker T_U$, $1\le j\le l$. It is also easy to see that the functions
$$
\U^*F^{-1}\bs{c}_1,\cdots,\U^*F^{-1}\bs{c}_r,\eta_1,\cdots,\eta_l
$$
span ${\frak S}_\Phi^{(\s_\vk)}$.

Note that the above reasoning is similar to the proof of Lemma 1.2 in \cite{PY2}.

Suppose now that $\Phi$ satisfies (C4). Let us show that it is very badly approximable.
As we have already observed (see Remark 2 after the statement of Theorem \ref{plo}),
the singular values $s_j(\Phi(\z))$ are constant almost everywhere on $\T$.
Let $\s_0>\cdots>\s_{\iota-1}$ 
be positive numbers such that for almost all $\z\in\T$ the distinct nonzero singular values of $\Phi(\z)$ 
are precisely $\s_0,\cdots,\s_{\iota-1}$.
We argue by induction on $\iota$. If $\iota=0$, the situation is trivial. Suppose that $\iota>0$.
Suppose that $\xi\in\Ker T_\Phi$ and $\xi(\z)$ is a maximizing vector of $\Phi(\z)$ for almost all $\z\in\T$.
Clearly, $H_\Phi\xi=\Phi\xi$ and $\|H_\Phi\xi\|=\s_0\|\xi\|$. It follows that $\xi$ is a maximizing vector 
of $H_\Phi$ and $\Phi$ is badly approximable. Conversely, if $\xi$ is a maximizing vector of $H_\Phi$, then
$\xi\in\Ker T_\Phi$ and $\xi(\z)$ is a maximizing vector of $\Phi(\z)$ for almost all $\z\in\T$
(see \S 2).

Suppose that $\dim{\frak S}_\Phi^{(\s_0)}(\z)=r$ for almost all $\z\in\T$. Obviously,
$\dim{\frak S}_{\Phi^{\rm t}}^{(\s_0)}(\z)=r$ for almost all $\z\in\T$. 
Let us first show that $\Phi$ admits a partial canonical factorization \rf{pcf} in which $\V$ and $\W^{\rm t}$ 
are $r$-balanced unitary-valued matrix functions, $U$ is an $r\times r$ 
very badly approximable unitary-valued matrix function such that $\|H_U\|_{\rm e}<1$.

It is well known and it is easy to verify that if $\xi$ is a maximizing vector of $H_\Phi$ and
$\eta=\bar z\ov{H_\Phi\xi}$, then $\eta$ is a maximizing vector of $H_{\Phi^{\rm t}}$ and vice versa.

Let $\cL$ be the minimal invariant subspace of multiplication by $z$ on $H^2(\C^n)$
that contains all maximizing vectors of $H_\Phi$ and let $\M$ be the minimal invariant subspace of multiplication 
by $z$ on $H^2(\C^m)$ that contains all maximizing vectors of $H_{\Phi^{\rm t}}$. 

By Theorem B, there exist $n\times r$
inner and co-outer matrix
functions $\U$ and $\O$ such that $\cL=\U H^2(\C^r)$ and $\M=\O H^2(\C^r)$.
Let $\Theta$ and $\Xi$ are inner and co-outer matrix 
functions such that the matrix functions
$$
\V=\left(\begin{array}{cc}\U&\ov{\Theta}\end{array}\right)\quad\mbox{and}\quad
\W^{\rm t}=\left(\begin{array}{cc}\O&\ov{\Xi}\end{array}\right)
$$
are $r$-balanced. Then $\Phi$ admits a factorization \rf{pcf} (see Remark 3 after Theorem D in \S 2).
Moreover, to show that $\Phi$ is very badly approximable, it suffices to verify that $\Psi$ is very badly
approximable. Clearly, $\Psi$ satisfies (C1). Let us verify that $\Psi$ satisfies (C4).

Clearly, for almost all $\z\in\T$, $\s_1,\cdots,\s_{\iota-1}$ are all distinct nonzero singular values of 
$\Psi(\z)$. If $\iota=1$, then $\Psi=\bs{0}$, and so $\Phi$ is very badly approximable (see Theorem D).
Suppose now that $\iota>1$. Consider the function ${\frak S}_\Psi^{(\s_d)}$, $1\le d\le\iota-1$. 

Let $\xi_1,\cdots,\xi_\vk$ be functions in $\Ker T_\Phi$ 
such that 
$$
{\frak S}_\Phi^{(\s_d)}(\z)=\spn\{\xi_1(\z),\cdots,\xi_\vk(\z)\}\quad\mbox{for almost all}\quad\z\in\T.
$$
Since $\V$ is unitary-valued, we have
$$
{\frak S}_{\Phi\V}^{(\s_d)}(\z)=\spn\{(\V^*\xi_1)(\z),\cdots,(\V^*\xi_\kappa)(\z)\}
\quad\mbox{for almost all}\quad\z\in\T.
$$
We have $\V^*\xi_j=\left(\begin{array}{c}\U^*\xi_j\\\Theta^{\rm t}\xi_j\end{array}\right)$. Obviously,
$\left(\begin{array}{c}(\U^*\xi_j)(\z)\\0\end{array}\right)$ is a maximizing vector of
$(\Phi\V)(\z)$, and so it belongs to ${\frak S}_{\Phi\V}^{(\s_d)}(\z)$. Thus 
$\left(\begin{array}{c}0\\(\Theta^{\rm t}\xi_j)(\z)\end{array}\right)\in{\frak S}_{\Phi\V}^{(\s_d)}(\z)$. 
It is easy to see that $(\Theta^{\rm t}\xi_j)(\z)\in{\frak S}_{\Psi}^{(\s_d)}(\z)$. 
Moreover, it is evident that
$$
{\frak S}_{\Psi}^{(\s_d)}(\z)=\spn\{(\Theta^{\rm t}\xi_j)(\z):~1\le j\le\vk\}\quad
\mbox{for almost all}\quad\z\in\T.
$$
Let us show that $\Theta^{\rm t}\xi_j\in\Ker T_\Psi$. 

It follows from \rf{pcf} that 
$$
\Psi=\Xi^*\Phi\ov{\Theta}.
$$
Hence,
$$
\Psi\Theta^{\rm t}\xi_j=\Xi^*\Phi\ov{\Theta}\Theta^{\rm t}\xi_j=\Xi^*\Phi\xi_j\in H^2_-(\C^m),
$$
since $\xi_j\in\Ker T_\Phi$. By the inductive hypothesis, $\Psi$ is very badly approximable, and so
$\Phi$ is very badly approximable. $\bl$

\

\section{\bf An alternative approach: weighted estimates \\ and superoptimal weights}
\setcounter{equation}{0}

\

In this section we present an alternative, more geometric  proof of  the main result (Theorem
\ref{plo}).  Main ideas  of this proof go back to \cite{T}, where the so-called
 {\em superoptimal weights} were used to prove the uniqueness of superoptimal approximation. 

Although we do not use superoptimal weight per se in this proof, the main ideas from \cite{T}
(weighted estimates, optimal vectors, ``pinching'' the weights, etc) are present here, so we wanted
to mention the origin of the ideas. 

\medskip

{\bf 5.1. Matrix  weights and weighted Nehari Problem.} 
Let $W$ be  an $n\times n$ {\em matrix weight}, i.e., a bounded matrix-valued function on
$\T$, whose  values are nonnegative $n\times n$ matrices. Given a matrix weight,
one can introduce the weighted norm $\|\cdot\|_W$ on $L^2(\C^n)$:
$$
\|f\|_W^2 := (Wf, f)_{L^2(\C^n)} = \int_\T (W(\xi)f(\xi), f(\xi))dm(\xi), \qquad f\in
L^2(\C^n),
$$
with the corresponding weighted inner product 
$(\,\cdot\,, \,\cdot\,)_W$, $(f, g)_W = (W f,g)_{L^2(\C^n)}$. 

Given a Hankel
operator
$H_\Phi :H^2(\C^n)\to H^2(\C^m)$, we call the weight $W$  admissible (for the Hankel operator
$H_\Phi$) if the following inequality
$$
\|H_\Phi f\|^2\le \|f\|_W^2 = (Wf,f)_{L^2(\C^n)}:=
\int_{\T}(W(\xi)f(\xi),f(\xi))dm(\xi),\qquad 
f\in H^2(\C^n),
$$
holds. 

We need the following weighted analogue of the classical Nehari
Theorem.

\medskip

{\bf Weighted Nehari Theorem.}
{\it Let $\Phi\in L^\infty(\mm_{m,n})$ and let $W\in L^\infty(\mm_{n,n})$ be an
admissible weight for $\Phi$. Then there exists $F\in H^\infty(\mm_{m,n})$ such that
 the function $\Psi = \Phi - F$ satisfies the inequality
$
\Psi(\xi)^*\Psi(\xi) \le W(\xi)
$
a.e.~on $\T$.}

\medskip

This  theorem (and even its operator-valued version) easily follows from the classical
operator  Nehari Theorem. We refer the reader to \cite{T} for the proof.

\medskip

{\bf 5.2. The necessity of condition (C4).}  Suppose that
$\Psi$ is a very badly approximable function. By \rf{sv},
$s_k(\Phi(\z))=t_k(\Phi)$ for almost all $\z\in\T$, where the
$t_k(\Phi)$ are the {\em superoptimal singular values} of $\Phi$.  Let
$\s_k$, $k=0, 1, 2, \cdots$, be the sequence of {\em distinct}
nonzero superoptimal  singular values of $\Phi$ arranged in the decreasing
order. In other words, for almost all $\z\in\T$, the sequence $\s_k$,
$k=0, 1, 2, \cdots$, is the sequence of distinct singular values of
$\Phi(\z)$ arranged in the decreasing order. 

Define the functions $\f_k$ by $\f_k(x)= \max\{x, \s_k^2\}$, $x\ge0$,  
$k=0, 1, 2, \cdots$, and define the weights $W_k$ by  $W_k(\z)=\f_k(\Phi(\z)^*\Phi(\z))$, $\z\in\T$. 

Since $\Phi^*\Phi\le W_k$, the weights $W_k$ are admissible for the Hankel
operator $H_\Phi$. 

For $k=0,1,\cdots$, we denote by $\cE_k$ the set of all extremal functions
for the weighted estimate $\|H_\Phi f\|^2 \le (W_k f, f)$, i.e., the set of
all functions 
$f\in H^2(\C^n)$ satisfying
$$
\| H_\Phi f\|_2^2 = (W_kf, f)_{L^2(\C^n)}. 
$$

Since $\|H_\Phi\|=t_0(\Phi) > \|H_\Phi\|_{\rm e}$, the norm of $H_\Phi$ is
attained, and  $\cE_0$ is a nontrivial finite-dimensional subspace of
$H^2(\C^n)$. Since by the assumption of the theorem $\s_k>\|H_\Phi\|_{\rm
e}$, the subspaces $\cE_k$ are finite-dimensional, and since the sequence
$\cE_k$ is clearly increasing, all $\cE_k$ are nontrivial subspaces.

Denote by $E_k(\z)\df \spn\{f(\z):f\in\E_k\}$, $\z\in \T$.  More 
precisely, take some basis in $\E_k$, select a function $f_j$ from each
equivalence class, and define \lb$E_k(\z)=\spn\{f_j(\z):j=1,2,\cdots\}$. Note
that different choices of bases and representatives  give us different
functions $E_k$, but any two such functions coincide almost everywhere.
Thus the corresponding equivalence class of subspace-valued functions is
well defined.

It is easy to show that the function $\dim E_k(z)$ is constant almost everywhere on $\T$
and that the projection-valued functions $\z\mapsto P_{E_k(\z)}$ are measurable, cf \cite{T}.  

Our goal is to show, that  $E_k(\z)=  {\frak S}_\Phi^{(\s_k)}(\z)$ for almost all
$\z\in\T$. Then we are done, because any $\cE_k\subset\Ker T_\Phi$.
Indeed, for $f\in \cE_k$,
$$
(\Phi^*\Phi f, f)\le (W_k f, f) = \|H_\Phi f\|^2 \le (\Phi^*\Phi f, f), 
$$
whence $\|\Phi f\| = \|H_\Phi f\|$. Keeping in mind that $\|\Phi f\|^2 =
\|H_\Phi f\|^2 + \|T_\Phi f\|^2$, we get $\|T_\Phi f \|=0$.

Let us show first that $E_k(\z)\subset  {\frak S}_\Phi^{(\s_k)}(\z)$
for almost all $\z\in\T$. Assume the contrary. Then there exists a function $f\in \cE_k$
such that $f(\z)\notin {\frak S}_\Phi^{(\s_k)}(\z)$ on a set of
positive measure. Since for any {\em finite} collection of
functions  $f_1, f_2, \cdots, f_N\in H^2(\C^n)$,\lb the dimension $\dim\spn\{f_1,
f_2, \cdots, f_N\}$ is constant almost everywhere on $\T$ (the minors belong to the
Nevanlinna class), it follows that 
$f(\z)\notin {\frak S}_\Phi^{(\s_k)}(\z)$ for almost all $\z\in\T$.
By the definition of ${\frak S}_\Phi^{(\s_k)}$, we have 
$\|\Phi(\z) f(\z)\|_{\C^n} < \s_k \|f(\z)\|_{\C^n}$, $\z\in\T$, and so 
$$
\|H_\Phi f\|_2 \le \|\Phi f\|_2 < \s_k \|f\|_2. 
$$
However, this contradicts the definition of $\cE_k$. Hence, 
$E_k(\z)\subset  {\frak S}_\Phi^{(\s_k)}(\z)$ for almost all $\z\in\T$. 

Let us now prove that  $ E_k(\z) = {\frak S}_\Phi^{(\s_k)}(\z)$. Suppose that 
$E_k(\z)$ is a proper subspace of  ${\frak S}_\Phi^{(\s_k)}(\z)$ for almost all
$\z\in\T$. Let us show that in this case $\Phi$ is not a very badly approximable
function. 

Let $N$ be the largest integer such that $s_N(\Phi(\z)) =\s_k$
for almost all $\z\in\T$ (recall that the functions $\z\mapsto s_j(\Phi(\z))$ are constant 
almost everywhere). This
means that for almost all $\z$ there are exactly $N+1$ singular values of $\Phi(\z)$ (counting
multiplicities) that are greater than or equal to $\s_k$. 

We want to construct a function $\Psi$ such that $\Phi -\Psi\in H^\infty(\mm_{m,n})$,
$$
s_j(\Psi(\z)) \le s_j(\Phi(\z)), \ \ \z\in\T \ \ \text{ for }j<N,
$$
but 
$$
s_N(\Psi(\z)) < s_N(\Phi(\z)),\quad\z\in\T,
$$ 
which would imply that $\Phi$ is not very badly approximable.

To do that we ``pinch'' the weight $W_k(\xi)$ in the directions orthogonal to $E_k(\xi)$ to
make it smaller (but still admissible)  and then solve the Weighted Nehari Problem.

Namely, consider the one-parametric family of weights
$W_k^{[a]}$, $a>0$, defined by
$$
W_k^{[a]}(\xi)= P_{E_k(\z)} W_k P_{E_k(\z)} \oplus a^2
P_{E_k(\z)^\perp} ;
$$
here we use the symbol $\oplus$ to emphasize that both operators on the right-hand side act on
orthogonal subspaces of $\C^n$, i.e., the operators $W_k^{[a]}(\xi)$ have
block-diagonal form with respect to  the orthogonal decomposition
$\C^n=E_k(\z) \oplus E_k(\z)^\perp$. 

If we can show that for some $a<\s_k$ the weight $W_k^{[a]}$ is still admissible, the
necessity is proved. Indeed, let $\Psi$ be a solution of the Weighted Nehari Problem, i.e., a
function such that $\Phi- \Psi \in H^\infty$ and $\Psi^*(\z)\Psi(\z)\le W_k^{[a]}(\z)$,
$\z\T$. Then the minimax property of the singular values implies that for $j<N$ 
$$
s_j(\Psi(\z)) \le s_j ( W_k^{[a]}(\z))^{1/2} \le s_j ( W_k(\z))^{1/2} =
s_j(\Phi(\z)),\quad\z\in\T,
$$
but 
$$
s_N(\Psi(\z))\le a < s_N(\Phi(\z)),\quad\z\in\T.
$$

We will need the following simple fact, whose proof is left as an exercise.  

\begin{lem}
\label{l6.3}
Let $T$ be an operator  (acting from one Hilbert space to
another one), and  let $f$ be a maximizing vector of $T$. Then for any vector $g$, the condition $g\perp f$
implies
$Tg\perp Tf$. 
\end{lem}

Let us now apply this lemma. We treat the Hankel operator $H_\Phi$ as
a operator, acting from $H^2(\C^n)$ endowed with the weighted norm
$\|\cdot\|_{W_k}$ to the space $H^2_-(\C^m)$. The nonzero vectors in $\cE_k$ are exactly
the maximizing vectors for this operator. Therefore by Lemma
\ref{l6.3}, for  any function $g\in H^2(\C^n)$ orthogonal to $\cE_k$ with respect to the
weighted inner product $(\cdot, \cdot)_W$, we have 
$$
H_\Phi g \perp H_\Phi \E_k
$$
(with respect to the usual, unweighted scalar product).

Put
\begin{equation}
\label{6.1}
q=\sup\left\{\|H_\Phi f\|:~\|f\|_{W_k}=1,~\mbox{and}~ \ f \text{ is $W_k$-orthogonal to }
\cE_k\right\}. 
\end{equation}
Since $W_k$ is an admissible weight, $q\le 1$. Moreover, the following lemma says that actually
$q<1$. 

\begin{lem}
\label{l6.4}
Let $W$ be an invertible admissible weight for a Hankel operator $H_\Phi$ such that
$W(\xi)
\ge a^2I$, $a>\|H_\Phi\|_{\rm e}$, and let $K$ be a closed subspace of $H^2(\C^n)$. If
$$
q=\sup\{ \| H_\Phi f\| : f\in K, \|f\|_W \le 1\} = 1, \footnote{Note that the supremum is 
always at most 1}
$$
then there exists a vector $f_0\in K$ such that $\| H_\Phi f_0\| = \|f_0\|_W$.
\end{lem}

\Pf
Putting $g= W^{1/2} f$, we can rewrite the  condition $q=1$ in the following way: 
$$
\sup\{ \|H_\Phi W^{-1/2} g :~ g\in W^{1/2} K, \|g\|=1\} =1, 
$$ 
which simply means that the norm of the operator $(H_\Phi W^{-1/2}) \big| W^{1/2}K$ is $1$. 
Since the norm of multiplication by $W^{-1/2}$ is at most $a^{-1}$, the essential norm of the operator
$(H_\Phi W^{-1/2}) \big| W^{1/2}K$ can be estimated as 
$$
\left\| (H_\Phi W^{-1/2}) \big| W^{1/2}K\right\|_{\rm e} \le \|H_\Phi\|_{\rm e} \|W^{-1/2}\|_\infty \le 
\|H_\Phi\|_{\rm e} a^{-1} <1.
$$
Therefore the norm of this operator is attained on some vector $g_0\in W^{1/2}K$, and so $f_0 = W^{-1/2}
g_0$ is a maximizing vector in $K$. $\bl$

Let us apply Lemma \ref{l6.4} to the weight $W_k$ and the subspace $K$ of $H^2(\C^n)$
of vectors that are $W_k$-orthogonal to $\cE_k$. If $q=1$ in \eqref{6.1}, the
lemma asserts that there is a maximizing vector in $K$, which is impossible, since $\cE_k$
contains {\em all} maximizing vectors. 

To complete the proof of necessity, we have to show that for $a=qa_k$, the weight $W_k$ is
still admissible. First of all, note that  $E_k(\z)$ is an invariant subspace of
all $W_k^{[a]}(\z)$ (including $W_k(\z)$) for almost all  $\z\in\T$. Since  for any $f\in \cE_k$,
we have $f(\z)\in E_k(\z)$, $\z\in\T$, and so for $f\in \cE_k$ and $g\in H^2(\C^n)$, we obtain 
\begin{multline}
\label{6.2}
(W_k^{[a]}(\z) f(\z), g(\z)) = (P_{E_k(\z)} W_k^{[a]}(\z) f(\z), g(\z)) \\
=(P_{E_k(\z)} W_k(\z) f(\z), g(\z)) = ( W_k (\z) f(\z), g(\z)),\quad\z\in\T. 
\end{multline}
In particular, it follows that $K$, being the $W_k$-orthogonal complement of $\cE_k$, is also
the $W_k^{[a]}$-orthogonal complement of $\cE_k$ for all $a>0$. 

Let $f\in \cE_k$ and let $g$ be $W_k$ orthogonal to $\cE_k$. Then $\|f\|_{W_k^{[a]}} $ does
not depend on $a$, and for $a=q a_k$ we have $q\|g\|_{W_k} \le \|g \|_{W_k^{[a]}}$. (If
$g(\z)$ were pointwise orthogonal to $E_k(\z)$, then equality would hold. But $g(\z)$ is not
necessarily pointwise orthogonal to $E_k(\z)$, so we can guarantee only inequality). By
Lemma \ref{l6.3}, $H_\Phi f \perp H_\Phi g$, and so 
\begin{multline*}
\| H_\Phi (f+g)\|^2 =\| H_\Phi f\|^2 + \| H_\Phi g\|^2 
\\  \le \|f\|^2_{W_k} + q\|g\|^2_{W_k}
\le 
\|f\|^2_{W_k^{[a]}} + \|g\|^2_{W_k^{[a]}} = \|f+g\|^2_{W_k^{[a]}}, 
\end{multline*}
whence the weight $W_k^{[a]}$ is admissible. This completes the proof of necessity. $\bl$

\medskip

{\bf 5.2. Sufficiency.} Suppose that a function $\Phi$ satisfies condition (C4).
Let us show that $\Phi$ is very badly approximable. 

As we already discussed above, (C4) implies that singular values of
$\Phi(z)$ (i.e.~the functions $\z \mapsto s_k(\Phi(\z))$) are constant almost everywhere on $\T$.
Let $s_0, s_1, s_2, \cdots$ denote these singular values arranged in the nonincreasing
order (counting multiplicity), and let $\s_0, \s_1, \s_2, \cdots$ be all {\em distinct} 
singular values  arranged in the decreasing order (i.e., $\s_0, \s_1, \s_2, \cdots$ be the
singular values of $\Phi(\z)$ {\em not counting multiplicity}).

Let $F$ be a superoptimal approximation of $\Phi$, $\Psi= \Phi- F$, and let $t_0,
t_1, t_2, \cdots$ be the superoptimal singular values of $\Phi$ (equivalently, of $\Psi$).
Let
$N_k$ be the largest integer such that  $s_{N_k} =\sigma_k$, which means that there are
exactly $N_k+1$ singular values that are greater than or equal to $\sigma_k$.   

As in the proof of necessity, let us introduce the weight $W=\Phi^*\Phi$, and let
$W_k(\xi)= \f_k(W)$, where $\f_k(x) = \max \{x, \s_k\}$ for $x\ge0$. 

We are going to prove using induction on $k$ that for all $k$ the following conditions are satisfied:

(i) $\Psi^*(\z)\Psi(\z) \le W_k(\xi)$ for almost all $\z\in\T$;

(ii) $\Psi_k(\z) \big|  {\frak S}_\Phi^{(\s_k)}(\z) = \Phi_k(\z) \big|  {\frak
S}_\Phi^{(\s_k)}(\z)$ for almost all $\z\in \T$;

(iii) $t_j =s_j$ for $0\le j\le N_k$.

This will immediately prove that $\Psi\equiv \Phi$, and so $\0$ is the unique superoptimal
approximation of $\Phi$. 

Consider first the case $k=0$. By the definition of {\em superoptimal approximation}
$t_0\le s_0=\s_0$, and hence,
$$
\Psi^*(\z)\Psi(\z) \le s_0^2 I = W_0(\z),\quad \z\in\T, 
$$
i.e., condition (i) is satisfied. 

Suppose that  $f_0, f_1, \cdots, f_{N_0}\in \Ker T_\Phi$ are functions that span ${\frak
S}_\Phi^{(\s_0)}$. Since
\begin{equation}
\label{6.3}
\Phi f = T_\Phi f + H_\Phi f, 
\end{equation}
and $T_\Phi f_j=0$, we have $H_\Phi f_j = \Phi f_j$ for $0\le j \le N_0$. Hence,
$$
\s_0 \|f_j\| = \|\Phi f_j\| = \|H_\Phi f_j \| = \|H_\Psi f_j \| \le \|\Psi f_j \| \le
t_0 \|f_j\|.
$$
Since $t_0\le \s_0$, the above inequalities are actually equalities
and \eqref{6.3} implies that 
$$
\Psi f_j = H_\Psi f_j = H_\Phi f_j = \Phi f_j. 
$$
Since $\spn\{f_j(\z): 0\le j \le N_0 \} = {\frak S}_\Phi^{(\s_0)}(\z)$ for almost all $\z\in\T$, 
condition (ii) is satisfied. Condition (iii) is an immediate consequence of (i) and (ii).

Let us assume now that the inductive hypotheses (i)--(iii) are proved for $k$, and we want to
prove them for $k+1$. It follows from (iii) and the definition of superoptimal approximation
that $t_{N_k+1} \le s_{N_k+1} = \s_{k+1}$, and so $\Psi^*\Psi \le W_{k+1}$.
This proves (i).

The proof of the other two condition is very similar to that of in the case $k=0$. 

First of all note that the case $\s_{k+1}=0$ is trivial, since in this case
$W_{k+1}(\z)$, $\Phi(\z)$ and $\Psi(\z)$ must be zero on ${\frak
S}_\Phi^{(\s_k)}(\z)^\perp$.

Let us assume
that $\s_{k+1}>0$ and let $f_j$,
$0\le j
\le N_{k+1}$ be functions in $\Ker T_\Phi$ that span 
${\frak S}_\Phi^{(\s_{k+1})}$. The condition $T_\Phi f_j =0$ and \eqref{6.3} implies that 
$H_\Phi f_j = \Phi f_j$ and using the fact that $f_j(\z) \in {\frak
S}_\Phi^{(\s_{k+1})}(\z)$ almost everywhere on $\T$, we can write
$$
(W_{k+1} f_j, f_j) = \|  \Phi f_j \|^2 = \|H_\Phi f_j\|^2 =  \|H_\Psi
f_j\|^2
\le
 \|\Psi f_j\|^2 = (\Psi^*\Psi f_j, f_j). 
$$
We have already proved that $\Psi^*\Psi \le W_{k+1}$, and so the inequality in the above chain
turns into equality. Thus \eqref{6.3} implies that
$$
\Psi f_j = H_\Psi f_j = H_\Phi f_j = \Phi f_j,  
$$
which in turn implies condition (ii) follows, since 
$$
\spn\{ f_j(\z): 0\le j \le N_{k+1}\} = {\frak
S}_\Phi^{(\s_{k+1})}(\z),\quad\z\in\T.
$$ 
Condition (iii) is again an immediate
consequence of (i) and (ii).  $\bl$ 

\

\section{\bf Badly Approximable Matrix Functions}
\label{s5}
\setcounter{equation}{0}

\

In this section we obtain a characterization of the badly approximable matrix functions
$\Phi$ satisfying the condition $\|H_\Phi\|_{\rm e}<\|H_\Phi\|$. Finally,
under the same assumption we characterize 
matrix functions $\Phi$, for which $\0$ is the only best approximation.

\begin{thm}
\label{baf}
Let $\Phi$ be a matrix function in $L^\be(\mm_{m,n})$ such that \lb$\|H_\Phi\|_{\rm e}<\|\Phi\|_{L^\be}$.
Then $\Phi$ is badly approximable if and only if the following conditions are satisfied:
\item{\rm(i)}
$\|\Phi(\z)\|_{\mm_{m,n}}$ is constant for almost all $\z\in\T$;
\item{\rm(ii)} 
there exists a function $f$ in $\Ker T_\Phi$ such that $f(\z)$ is
a maximizing vector of $\Phi(\z)$ for almost all $\z\in\T$.
\end{thm}

{\bf Remark 1.} It will be clear from the proof that if $\Phi$ is an arbitrary matrix function
satisfying (i) and (ii), then it is badly approximable. In other words, to prove that (i) and (ii) imply
that $\Phi$ is badly approximable, we do not need the condition $\|H_\Phi\|_{\rm e}<\|\Phi\|_{L^\be}$.

\medskip

\Pf Suppose that $\Phi$ is badly approximable. Then it admits a factorization
$$
\Phi=W^*\left(\begin{array}{cc}\s u&\bs{0}\\\bs{0}&\Psi\end{array}\right)V^*
$$
where $V$ and $W^{\rm t}$ are thematic (1-balanced) matrix functions, $\s=\|\Phi\|_{L^\be}$,
$u$ is a scalar unimodular badly approximable function such that $\|H_u\|_{\rm e}<1$, and
$\Psi$ is an $(m-1)\times(n-1)$ matrix function such that  $\|\Psi\|_{L^\be}\le\s$
(see \cite {AP}or \cite{P}, Ch. 14, \S 4). Let
$$
V=\left(\begin{array}{cc}\bs{v}&\ov{\Theta}\end{array}\right),\quad
W=\left(\begin{array}{cc}\bs{w}&\ov{\Xi}\end{array}\right)^{\rm t},
$$
where $\bs{v}$ and $\bs{w}$ are inner and co-outer column functions while $\Theta$ and $\Xi$ are 
inner and co-outer matrix functions. 

It follows from the characterization of badly approximable scalar functions \lb mentioned in the introduction
that $T_u$ is Fredholm and $\ind T_u>0$. Therefore \lb$\Ker T_u\ne\{\bs{0}\}$. Let $h$ be a nonzero function in $\Ker T_u$.
Put $f=h\bs{v}$. We have
\bey
\Phi f&=&W^*\left(\begin{array}{cc}\s u&\bs{0}\\\bs{0}&\Psi\end{array}\right)
\left(\begin{array}{c}\bs{v}^*\\\Theta^{\rm t}\end{array}\right)h\bs{v}
=W^*\left(\begin{array}{cc}\s u&\bs{0}\\\bs{0}&\Psi\end{array}\right)
\left(\begin{array}{c}h\\\bs{0}\end{array}\right)\\[.2cm]
&=&\left(\begin{array}{cc}\ov{\bs{w}}&\Xi\end{array}\right)
\left(\begin{array}{c}\s uh\\\bs{0}\end{array}\right)=
\s uh\bs{w}\in H^2_-(\C^m),
\eey
since $h\in\Ker T_u$. Thus $f\in\Ker T_\Phi$. On the other hand,
$$
\|\Phi(\z)f(\z)\|_{\C^m}=|\s u(\z)h(\z)|=\s|h(\z)|=\s\|f(\z)\|_{\C^n}
$$
for almost all $\z\in\T$, i.e., $f(\z)$ is a maximizing vector of $\Phi(\z)$ for almost all $\z\in\T$.

Suppose now that (i) holds, $f\in\Ker T_\Phi$, and $f(\z)$ is a maximizing vector of $\Phi(\z)$ for almost all $\z\in\T$.
Then $H_\Phi f=\Phi f$, and it is easy to see that 
$$
\|H_\Phi f\|=\|\Phi f\|=\|\Phi\|_{L^\be}\|f\|,
$$
i.e., $\|H_\Phi\|=\|\Phi\|_{L^\be}$, and so $\Phi$ is badly approximable. $\bl$

The following theorem describes {\em badly approximable} functions, for which $\0$ is the only best
approximation. If $\Phi$ is a nonzero matrix function in $L^\be(\mm_{m,n})$, we can normalize it by the
condition $\|\Phi\|_{L^\be}=1$.

\begin{thm}
\label{t5.2}
Let $\Phi$ be a matrix function in $L^\be(\mm_{m,n})$ such that \lb$\|H_\Phi\|_{\rm e}<\|\Phi\|_{L^\be}$.
Then $\0$ is the only best approximation of $\Phi$ if and only if the following conditions are satisfied:
\begin{enumerate}
\item[\rm(i)] $\Phi$ takes isometric values if $n\le m$ and $\Phi^{\rm t}$ takes isometric values if 
 $n>m$ almost everywhere on $\T$; 

\item[\rm(ii)] the function $\z\mapsto(\Ker \Phi(\z))^\perp$, $\z\in \T$, is 
spanned by finitely many functions in $\Ker T_\Phi$.
\end{enumerate}
\end{thm}

{\bf Remark 2.} If $n\le m$ and $\Phi$ satisfies (i), then $\Ker\Phi(\z)=\{0\}$ for almost all $\z\in\T$,
and so  (ii) means that there are finitely many functions $f_j\in\Ker T_\Phi$ such that
$\spn\{f_j(\z):~j=1,2,\cdots\}=\C^n$ for almost all $\z\in\T$. Note that if $n>m$, then instead of $\Phi$ we
can consider the transposed function $\Phi^{\rm t}$.

\medskip
{\bf Remark 3.} As in the case of Theorem \ref{baf}, to prove that (i) and (ii) imply that $\0$ is the only
best 
approximation of $\Phi$, we do not need the condition $\|H_\Phi\|_{\rm e}<\|\Phi\|_{L^\be}$.

\medskip

{\bf Proof.} Let $\Phi$ be a badly approximable matrix function in
$L^\be(\mm_{m,n})$ such that \lb$\|H_\Phi\|_{\rm e}<\|\Phi\|_{L^\be}$=1.
Let $r$ be the number of superoptimal
singular values of $\Phi$ equal to $\|\Phi\|_{L^\be}=1$. Suppose that $r<\min\{m,n\}$.
By Theorem D in \S \ref{s2}, $\Phi$ admits a factorization
$$
\Phi=\W^*\left(\begin{array}{cc} U&\bs{0}\\\bs{0}&\Psi\end{array}\right)\V^*,
$$
where $\V$ and $\W^{\rm t}$ are $r$-balanced matrix functions, 
$U$ is an $r\times r$ very badly approximable unitary-valued function such 
that $\|H_U\|_{\rm e}<1$, and
$\Psi$ is an $(m-r)\times(n-r)$ matrix function such that 
$\|\Psi\|_{L^\be}\le1$, $\|H_\Psi\|<1$, and $\|H_\Psi\|_{\rm e}\le\|H_\Phi\|_{\rm e}$.

Since $\|H_\Psi\|<1$, there exist infinitely many matrix functions $F\in
H^\infty(\mm_{m-r,n-r})$ such that $\|\Psi -F\|_\infty \le 1$. Note, that if $F\ne\0$, then the
function
$$
\W^*\left(\begin{array}{cc} \bs{0}&\bs{0}\\\bs{0}&F\end{array}\right)
\V^*
$$
is a nonzero  function in $H^\infty(\mm_{m,n})$. Hence, $\Phi$ has infinitely many
best approximations. Thus
$r=\min\{m,n\}$, which means that (i) holds.

Note that any {\em superoptimal} approximation is a best approximation.
Thus if $\bs{0}$ is the only best approximation, it is also the only
superoptimal approximation. So $\Phi$ is a very badly approximable
function, and condition (ii) follows from Theorem \ref{plo}. 

Suppose now that a function $\Phi$ satisfies (i) and
(ii). Let $F$ be a best approximation of $\Phi$, and let $\Psi = \Phi
-F$. Let $f_j$ be functions in $\Ker T_\Phi$ that span the function
$(\Ker \Phi(z))^\perp$. The condition $T_\Phi f_j
=\bs{0}$ implies that $H_\Phi f_j = \Phi f_j$ (see \eqref{6.3}), and
so
$$
\|f_j\|_2 = \|\Phi f_j\|_2 =\| H_\Phi f_j \|_2 = \|H_\Psi f_j\|_2 \le \|
\Psi f_j \|_2 \le \|f_j\|_2. 
$$ 
Therefore all inequalities in the above chain must be equalities, and it
follows from \eqref{6.3} that 
$$
\Psi f_j = H_\Psi f_j = H_\Phi f_j = \Phi f_j.
$$
Hence,
$$
\Psi(\z) \big| (\Ker \Phi(\z))^\perp = \Psi(\z) \big| (\Ker
\Phi(\z))^\perp, \quad \z\in\T. 
$$
If $n\le m$, then $(\Ker \Phi(\z))^\perp = \C^n$, and therefore
$\Phi\equiv \Psi$. 

To show that $\Phi\equiv \Psi$ for $m<n$ one more step is needed. Namely,
let us observe that $\Psi(\z)$ are contractions and that
$\Phi(\z)$ are co-isometries (i.e., $\Phi(\z)^* $ are isometries) for almost all
$\z\in\T$. It follows from Lemma \ref{l6.3} that if a contraction $T$ and a
co-isometry $U$ coincide on $(\Ker U)^\perp$, then $T\big| \Ker U =\bs{0}$,
and so $T=U$. 

Thus we have proved that $\Psi = \Phi$, i.e., $F=\bs{0}$,  and so $\bs{0}$ is the
only best approximation of $\Phi$. $\bl$

\

\

\noindent
\begin{tabular}{p{8cm}p{14cm}}
V.V. Peller & S.R. Treil \\
Department of Mathematics & Department of Mathematics \\
Michigan State University  & Brown University \\
East Lansing, Michigan 48824 & Providence, Rhode Island 02912\\
USA&USA%\\
\end{tabular}

\end{document}